\documentclass[11pt]{article} \textheight = 24truecm \textwidth
= 16truecm \hoffset = -2truecm \voffset = -2truecm
\usepackage[english]{babel}
\usepackage{amssymb,amsthm}
\usepackage{amsmath}

\begin{document}

\large \sloppy

\begin{center}
\textsc {On perfect metrizability of the functor of idempotent
probability measures}

\textsl{A.A.Zaitov, Kh.F.Kholturayev}
\end{center}

\begin{abstract}
In this paper we establish that the functor of idempotent
probability measures acting in the category of compacta and their
continuous mappings is perfect metrizable.
\end{abstract}

\textit{Keywords}: metric, idempotent probability measures.

\textit{2000 Mathematics Subject Classification}. Primary 54C65,
52A30; Secondary 28A33.\\

The notion of idempotent measure finds important applications in
different parts of mathematics, mathematical physics, economics,
mathematical biology and others. One can find a row of
applications of idempotent mathematics from [1].

Consider the set $\mathbb{R}$ of real numbers with two algebraic
operations: addition $\oplus$ and multiplication $\odot$ defined
as $u\oplus v=\textrm{max}\{u,v\}$ and $u\odot v=u+v$.
$\mathbb{R}$ forms semifield with respect to this operations and,
the unity $\mathbf{1}=0$ and zero $\mathbf{0}=-\infty$, i. e.

($i$) the addition $\oplus$ and the multiplication $\odot$ are
associative;

($ii$) the addition $\oplus$ is commutative;

($iii$) the multiplication $\odot$ is distributive with respect to
the addition $\oplus$;

($iv$) each nonzero element $x\in \mathbb{R}$ is invertible.

It denotes by $\mathbb{R}_{\textrm{max}}$. It is idempotent, i. e.
$x\oplus x=x$ for all $x\in\mathbb{R}$, and commutative, i. e. the
multiplication $\odot$ is commutative.

Let $X$ be a compact Hausdorff space, $C(X)$ the algebra of
continuous functions $\varphi : X \rightarrow \mathbb{R}$ with the
usual algebraic operations. On $C(X)$ the operations $\oplus$ and
$\odot$ define as follow:

$\varphi \oplus \psi=\textrm{max}\{\varphi, \psi \}$, where
$\varphi, \psi \in C(X)$,

$\varphi \odot \psi=\varphi + \psi$, where $\varphi$, $\psi \in
C(X)$,

$\lambda \odot \varphi=\varphi+\lambda_{X}$, where $\varphi\in
C(X)$, $\lambda\in \mathbb{R},$  and $\lambda_X$ is a constant
function.

Recall [2] that a functional $\mu : C(X)\rightarrow
\mathbb{R}(\subset \mathbb{R}_{\textrm{max}})$ is called to be an
idempotent probability measure on $X$, if:

1)  $\mu (\lambda_{X})=\lambda$ for each $\lambda \in \mathbb{R}$;

2)  $\mu (\lambda \odot \varphi)=\mu (\varphi)+\lambda$ for all
$\lambda \in \mathbb{R}$, $\varphi \in C(X)$;

3)  $\mu (\varphi \oplus \psi)=\mu(\varphi)\oplus \mu(\psi)$ for
every $\varphi$, $\psi\in C(X)$.

For a compact Hausdorff space $X$ a set of all idempotent
probability measures on $X$ denotes by $I(X)$.  Consider $I(X)$ as
a subspace of $\mathbb{R}^{C(X)}$. In the induced topology the
sets of the view
\begin{center} $\langle \mu; \varphi_1, \varphi_2, ...,
\varphi_k; \varepsilon \rangle=\{\nu\in I(X):
|\mu(\varphi_i)-\nu(\varphi_i)|<\varepsilon, i=1, ..., k\}$,
\end{center}
form a base of neighborhoods of the idempotent measure $\mu\in
I(X)$, where $\varphi_i\in C(X)$, $i=1, ..., k$, and $\varepsilon
>0$. The topology generated by this base coincide with pointwise
topology on $I(X)$. The topological space $I(X)$ is compact [2].
Given a map $f:X\rightarrow Y$ of compact Hausdorff spaces the map
$I(f):I(X)\rightarrow I(Y)$ defines by the formula
$I(f)(\mu)(\varphi)=\mu(\varphi\circ f)$, $\mu\in I(X)$, where
$\varphi\in C(Y)$. The construction $I$ is a covariant functor,
acting in the category of compact Hausdorff spaces and their
continuous mappings. Moreover, $I$ is uniform metizable functor
[3].

Since $I$ is normal functor then for an arbitrary idempotent
measure $\mu\in I(X)$ we may define the support of $\mu$:
$\mbox{supp}\mu=\bigcap\{A\subset X: \overline{A}=A,$ $\mu\in
I(A)\}$. For brevity, put $S\mu=\mbox{supp}\mu$. For a positive
integer $n$ put $I_n(X)=\{\mu\in I(X): |S\mu|\leq n\}$. Put
$I_{\omega}(X)=\bigcup\limits_{n=1}^\infty I_n(X)$. It is known
[2] that $I_{\omega}(X)$ is everywhere dense in $I(X)$. An
idempotent probability measure $\mu\in I_{\omega}(X)$ is named as
an idempotent probability measure with finite support. Note that
if $\mu$ is an idempotent probability measure with a finite
support $S\mu=\{x_1, x_2, ..., x_k \}$ then it represents in the
form
$$
\mu=\lambda_1\odot\delta_{x_1}\oplus\lambda_2\odot\delta_{x_2}\oplus
...\oplus\lambda_k\odot\delta_{x_k} \eqno(1)
$$
uniquely, where $\lambda_i\in \mathbb{R}_{\textrm{max}}$, $i=1,
..., k$, $\lambda_1\oplus\lambda_2\oplus
...\oplus\lambda_k=\textbf{1}$.

Let $\mu_1,\ \mu_2\in I_{\omega}(X)$. Then by (1) we have
$\mu_k=\bigoplus\limits_{i=1}^{n_k}\lambda
_{k}\odot\delta_{x_{ki}}$, $i=1,\ 2$. Put
$$\Lambda_{12}=\Lambda(\mu_1,\ \mu_2)=\{\xi\in I(X^2): I(\pi_i)(\xi)=\mu_i,\ i=1,2\},$$
where $\pi_i:X\times X\longrightarrow X$ is projection onto $i$-th
factor, $i=1,2$. By definition for each idempotent probability
measure $\xi\in\Lambda(\mu_1,\ \mu_2)$ we have
$\bigoplus\limits_{(x_{1j}, x_{2k})\in
S\xi}|\lambda_{2k}-\lambda_{1j}|\odot\rho(x_{1j}, x_{2k})<\infty$.
On the other hand as the set
$\{|\lambda_{2k}-\lambda_{1j}|\odot\rho(x_{1j}, x_{2k}): j=1,...,
n_1; k=1,..., n_2\} $ is finite there exists the number
$\min\limits_{\xi\in\Lambda_{12}}\{\bigoplus\limits_{(x_{1j},
x_{2k})\in S\xi} |\lambda_{2k}-\lambda_{1j}|\odot\rho(x_{1j},
x_{2k})\}.$ Put
\begin{center}
$H(\mu_1,\
\mu_2)=\min\limits_{\xi\in\Lambda_{12}}\{\bigoplus\limits_{(x_{1j},
x_{2k})\in S\xi} |\lambda_{2k}-\lambda_{1j}|\odot\rho(x_{1j},
x_{2k})\}.$
\end{center}

In [4] shown that the map $\rho_\omega:I_\omega(X)\times
I_\omega(X)\rightarrow \mathbb{R}$ defined as $\rho_\omega(\mu_1,
\mu_2)=\min\{\textrm{diam} X, \ H(\mu_1,\ \mu_2)\}$ is metric.
Moreover, the function $\rho_I :I(X)\times I(X)\rightarrow
\mathbb{R}$ which is an extension of $\rho_\omega$ onto completion
$I(X)$ of $I_\omega(X)$ is a metric as well.

We note this paper is continuation of papers [3-4]

Let now recall some notions. Put $I^0(X)=X$,
$I^n(X)=I(I^{n-1}(X))$, $n=1,\ 2,\ ...\ $. Consider two systems
$\eta$ and $\psi$. The system $\eta$ consists of all maps $\eta_X:
X\rightarrow I(X)$, $X\in Comp$, where $\eta_X$ defines as
$\eta_X(x)=\delta_x$, $x\in X$. Here $\delta_x$ is a Dirac measure
concentrated on $\{x\}$. The system $\psi$ consists of all
mappings $\psi_X: I^2(X)\rightarrow I(X)$, acting as the
following. Given $M\in I^2(X)$ put
$\psi_X(M)(\varphi)=M(\overline{\varphi})$, where for any function
$\varphi\in C(X)$ the function $\overline{\varphi}:
I(X)\rightarrow \mathbb{R}$ defines by the formula
$\overline{\varphi}(\mu)=\mu(\varphi)$. Fix compactum $X$ and for
a positive integer $n$ put $\psi_{n+1, n}=\psi_{I^{n-1}(X)}:
I^{n+1}(X)\rightarrow I^n(X)$ and $\eta_{n,\ n+1}=\eta_{I^n(X)}:
I^n(X)\rightarrow I^{n+1}(X)$. Note that $\psi_{n+1,
n}\circ\eta_{n,\ n+1}=Id_{I^n(X)}$ [2].

\textbf{Lemma 1.} \textsl{$\psi_{0,\ 1}:I^2(X)\rightarrow I(X)$ is
non-expanding map.}

\textsc{Proof}. It is enough to consider idempotent probability
measures with everywhere finite supports. Let $M_1$ and $M_2$ are
such measures from $I^2(X)$ and let $SM_1=\{\mu_{11},\ \mu_{12},\
...,\ \mu_{1n_1}\}$, $SM_2=\{\mu_{21},\ \mu_{22},\ ...,\
\mu_{2n_2}\}$, be the their supports, where $\mu_{lk}$ are
idempotent probability measures with finite supports, $k=1,\ ...,\
n_l$, $l=1,\ 2$. Assume $S\mu_{lk}=\{x_{k1}^l,\ ...,\
x_{kt_k}^l\}$. Then we have
\begin{center}
$M_l=m_{l1}\odot\delta_{\mu_{l1}}\oplus
m_{l2}\odot\delta_{\mu_{l2}}\oplus\ ...\ \oplus
m_{ln_l}\odot\delta_{\mu_{ln_l}}, $\\
$\mu_{lk}=\lambda_{k1}^l\odot\delta_{x_{k1}^l}\oplus
\lambda_{k2}^l\odot\delta_{x_{k2}^l}\oplus\ ...\ \oplus
\lambda_{kt_k}^l\odot\delta_{x_{kt_k}^l},$
\end{center}
where $t_k$ are positive integers, $k=1,\ 2,\ ...,\ n_l$, $l=1,\
2$. By definition for any $\varphi\in C(X)$ one has
\begin{center}
$\psi_X(M_l)(\varphi)=M_l(\overline{\varphi})=(m_{l1}\odot\delta_{\mu_{l1}}\oplus
m_{l2}\odot\delta_{\mu_{l2}}\oplus\ ...\ \oplus
m_{ln_1}\odot\delta_{\mu_{ln_l}})(\overline{\varphi})= $\\
$=m_{l1}\odot\overline{\varphi}(\mu_{l1})\oplus
m_{l2}\odot\overline{\varphi}(\mu_{l2})\oplus\ ...\ \oplus
m_{ln_l}\odot\overline{\varphi}(\mu_{ln_l})= $\\
$=m_{l1}\odot\mu_{l1}(\varphi)\oplus
m_{l2}\odot\mu_{l2}(\varphi)\oplus\ ...\ \oplus
m_{ln_l}\odot\mu_{ln_l}(\varphi)= $\\
$=m_{l1}\odot(\lambda_{11}^l\odot\delta_{x_{11}^l}\oplus
\lambda_{12}^l\odot\delta_{x_{12}^l}\oplus\ ...\ \oplus
\lambda_{1t_1}^l\odot\delta_{x_{1t_1}^l}(\varphi))\oplus $\\
$\oplus m_{l2}\odot(\lambda_{21}^l\odot\delta_{x_{21}^l}\oplus
\lambda_{22}^l\odot\delta_{x_{22}^l}\oplus\ ...\ \oplus
\lambda_{2t_2}^l\odot\delta_{x_{2t_2}^l}(\varphi))\oplus\ ...\
\oplus $\\
$\oplus
m_{ln_l}\odot(\lambda_{n_l1}^l\odot\delta_{x_{n_l1}^l}\oplus
\lambda_{n_l2}^l\odot\delta_{x_{n_l2}^l}\oplus\ ...\ \oplus
\lambda_{n_lt_{n_l}}^l\odot\delta_{x_{n_lt_{n_l}}^l}(\varphi))=$\\
$=\left(\bigoplus\limits_{s=1}^{t_1}(m_{l1}\odot\lambda_{1s}^l)\odot\delta_{x_{1s}^l}
\oplus\bigoplus\limits_{s=1}^{t_2}(m_{l2}\odot\lambda_{2s}^l)\odot\delta_{x_{2s}^l}
\oplus\bigoplus\limits_{s=1}^{t_{n_l}}(m_{ln_l}\odot\lambda_{n_ls}^l)\odot\delta_{x_{n_ls}^l}\right)
(\varphi)=$\\
$=\left(\bigoplus\limits_{k=1}^{n_l}\bigoplus\limits_{s=1}^{t_k}(m_{lk}\odot\lambda_{ks}^l)
\odot\delta_{x_{ks}^l}\right)(\varphi),$
\end{center}
i. e.
$$\psi_X(M_l)=\bigoplus\limits_{k=1}^{n_l}\bigoplus\limits_{s=1}^{t_k}(m_{lk}\odot\lambda_{ks}^l)
\odot\delta_{x_{ks}^l}\qquad l=1,\ 2.\eqno(2)
$$
From (2) immediately follows that
$$S\psi_X(M_l)=\{x_{11}^l,
x_{12}^l, ..., x_{1t_1}^l, x_{21}^l, x_{22}^l, ..., x_{2t_2}^l,
..., x_{n_l1}^l, x_{n_l2}^l, ..., x_{n_lt_{n_l}}^l\},\qquad l=1,\
2.$$

Suppose now $\Xi$ be an existing according to Lemma 3 [3]
idempotent probability measure from $\Lambda(M_1,\ M_2)$ such that
\begin{center}
$\rho_{I^2}(M_1,\ M_2)=\min\{\bigoplus\limits_{(\mu_{1i},
\mu_{2j})\in S\Xi}|m_{1i}-m_{2j}|\odot\rho_I(\mu_{1i}, \mu_{2j}),\
\mbox{diam}I(X)\}.$
\end{center}
Let $\xi_{ij}\in \Lambda(\mu_{1i},\ \mu_{2j})$ be idempotent
measures existing by Lemma 3 [3]. Then since
$\mbox{diam}I(X)=\mbox{diam}X$ we have
\begin{center}
$\rho_{I^2}(M_1, M_2)=\min\{\bigoplus\limits_{(\mu_{1i},
\mu_{2j})\in
S\Xi}(|m_{1i}-m_{2j}|\odot\bigoplus\limits_{(x_{ip}^1,
x_{jq}^2)\in S\xi_{ij}}|\lambda_{ip}^1-
\lambda_{jq}^2|\odot\rho(x_{ip}^1, x_{jq}^2)), \mbox{diam}X\}=$\\
$=\min\{\bigoplus\limits_{(\mu_{1i}, \mu_{2j})\in
S\Xi}(\bigoplus\limits_{(x_{ip}^1, x_{jq}^2)\in
S\xi_{ij}}|m_{1i}-m_{2j}|\odot|\lambda_{ip}^1-
\lambda_{jq}^2|\odot\rho(x_{ip}^1, x_{jq}^2)),
\mbox{diam}X\}\geq$\\
$\geq\min\{\bigoplus\limits_{(\mu_{1i}, \mu_{2j})\in
S\Xi}(\bigoplus\limits_{(x_{ip}^1, x_{jq}^2)\in
S\xi_{ij}}|m_{1i}\odot\lambda_{ip}^1-m_{2j}\odot
\lambda_{jq}^2|\odot\rho(x_{ip}^1, x_{jq}^2)),
\mbox{diam}X\}\geq$\\
$\geq\min\{\min\limits_{\xi\in \Lambda(\psi_X(M_1),
\psi_X(M_2))}\{\bigoplus\limits_{(x_{ip}^1, x_{jq}^2)\in
S\xi_{ij}}|m_{1i}\odot\lambda_{ip}^1-m_{2j}\odot
\lambda_{jq}^2|\odot\rho(x_{ip}^1, x_{jq}^2)\},
\mbox{diam}X\}=$\\
$=\rho_I(\psi_X(M_1), \psi_X(M_2)),$
\end{center}
i. e. $\rho_{I^2}(M_1, M_2)\geq\rho_I(\psi_X(M_1), \psi_X(M_2))$.
Lemma 1 is proved.

\textbf{Lemma 2.} \textsl{$\rho_I(\mu,\
\delta_{x_0})=\rho_{I^2}(\delta_{\delta_{x_0}},\ N)$ for each
$N\in \psi_X^{-1}(\mu)$.}

\textsc{Proof}. According to our metric it is enough to consider
idempotent probability measures with everywhere finite supports.
Fix an arbitrary point $x_0\in X$ and let $\mu\in I_\omega(X)$ be
an arbitrary measure. Then $\mu$ has decomposition of the form
(1). It is easy to see that measure
$\mu\otimes\delta_{x_0}=\lambda_1\odot\delta_{(x_1,
x_0)}\oplus\lambda_2\odot\delta_{(x_2,
x_0)}\oplus\lambda_k\odot\delta_{(x_k, x_0)}$ is an unique element
of the set $\Lambda(\mu, \delta_{x_0})$. Consequently we have
$$\rho_I(\mu, \delta_{x_0})=\min\{\mbox{diam}X,\
|\lambda_1|\odot\rho(x_1, x_0)\oplus|\lambda_2|\odot\rho(x_2,
x_0)\oplus|\lambda_k|\odot\rho(x_k, x_0)\}.$$

Let now $N\in \psi_X^{-1}(\mu)$ be a measure with everywhere
finite supports. Then we have
\begin{center}
$N=\alpha_1\odot\delta_{\nu_1}\oplus
\alpha_2\odot\delta_{\nu_2}\oplus...\oplus \alpha_s\odot\delta_{\nu_s},$\\
$\nu_i=\lambda_1^i\odot\delta_{x_1^i}\oplus\lambda_2^i\odot\delta_{x_2^i}\oplus
...\oplus\lambda_{k_i}^i\odot\delta_{x_{k_i}^i}$,
\end{center}
where $i=1,\ 2,\ ...,\ s$. Hence $S\psi_X(N)=\{x_{j}^i:\ j=1,\ 2,\
...,\ k_i;\ i=1,\ 2,\ ...,\ s\}$. On the other hand
$S\psi_X(N)=S\mu=\{x_1,\ x_2,\ ...,\ x_k\}$. Put $J_l=\{x_{j}^i\in
S\psi_X(N):\ x_{j}^i=x_l\}$, $l=1,\ 2,\ ...,\ k$. After making
slight modifications we can (2) rewrite in the view
$\psi_X(N)=\bigoplus\limits_{l=1}^k(\bigoplus\limits_{x_j^i\in
J_l}\alpha_i\odot\lambda_j^i)\odot\delta_{x_l}$. Since
$\psi_X(N)=\mu$ then for each $l=1,\ 2,\ ...,\ k$,  we have
$\lambda_l=\bigoplus\limits_{x_j^i\in
J_l}\alpha_i\odot\lambda_j^i$.

Let now we find the distance between $N$ and
$\delta_{\delta_{x_0}}$:
\begin{center}
$\rho_{I^2}(N,\ \delta_{\delta_{x_0}})=\min\{\mbox{diam}I(X),\
\bigoplus\limits_{i=1}^s(|\alpha_i|\odot\rho_I(\nu_i,\
\delta_{x_0}))\}=$\\
$=\min\{\mbox{diam}I(X),\ \bigoplus\limits_{i=1}^s(|\alpha_i|\odot
\min\{\mbox{diam}X,\
\bigoplus\limits_{j=1}^{k_i}|\lambda_j^i|\odot\rho(x_j^i,\
x_0)\})\}=$\\
$=\min\{\mbox{diam}X,\ \bigoplus\limits_{i=1}^s(|\alpha_i|\odot
\bigoplus\limits_{j=1}^{k_i}|\lambda_j^i|\odot\rho(x_j^i,\
x_0))\}=$\\
$=\min\{\mbox{diam}X,\
\bigoplus\limits_{i=1}^s\bigoplus\limits_{j=1}^{k_i}(|\alpha_i|\odot
|\lambda_j^i|\odot\rho(x_j^i,\
x_0))\}=$\\
$=\min\{\mbox{diam}X,\
\bigoplus\limits_{i=1}^s\bigoplus\limits_{j=1}^{k_i}(|\alpha_i\odot
\lambda_j^i|\odot\rho(x_j^i,\
x_0))\}=$\\
$=\min\{\mbox{diam}X,\
\bigoplus\limits_{l=1}^k\bigoplus\limits_{x_j^i\in
J_l}(|\alpha_i\odot \lambda_j^i|\odot\rho(x_l,\
x_0))\}=$\\
$=\min\{\mbox{diam}X,\
\bigoplus\limits_{l=1}^k|\lambda_l|\odot\rho(x_l,\
x_0))\}=\rho_I(\mu,\ \delta_{x_0})$.
\end{center}
Lemma 2 is proved.

\textbf{Lemma 3.} \textsl{If $\rho_I(\mu,\ \eta_{0,\ 1}(X))\geq
\varepsilon$ then $\rho_{I^2}(I(\eta_{0,\ 1})(\mu),\ \eta_{1,\
2}(I(X)))\geq \varepsilon$.}

\textsc{Proof}. As the above it is enough to consider idempotent
probability measures with everywhere finite supports. Without loss
of generality we may assume $\varepsilon\leq\mbox{diam}X$. In this
case $\rho_{I^n}=H_{I^n}$ for all positive integers $n$.

Let $\mu=\bigoplus\limits_{l=1}^n\lambda_l\odot\delta_{x_l}$ be an
arbitrary idempotent probability measure with finite support. For
any $\Phi\in C(I(X))$ we have
\begin{center}
$I(\eta_{0,\ 1})(\mu)(\Phi)=I(\eta_{0,\
1})(\bigoplus\limits_{l=1}^n\lambda_l\odot\delta_{x_l})(\Phi)=
(\bigoplus\limits_{l=1}^n\lambda_l\odot\delta_{x_l})(\Phi\circ\eta_{0,\
1})=$\\
$=\bigoplus\limits_{l=1}^n\lambda_l\odot\delta_{x_l}(\Phi\circ\eta_{0,\
1})=\bigoplus\limits_{l=1}^n\lambda_l\odot\Phi(\eta_{0,\
1}(x_l)=\bigoplus\limits_{l=1}^n\lambda_l\odot\Phi(\delta_{x_l})
=(\bigoplus\limits_{l=1}^n\lambda_l\odot\delta_{\delta_{x_l}})(\Phi)$.
\end{center}
This means that $I(\eta_{0,\ 1})(\mu)=
\bigoplus\limits_{l=1}^n\lambda_l\odot\delta_{\delta_{x_l}}$.

Let now $N=\eta_{1,\
2}(\bigoplus\limits_{i=1}^s\alpha_i\odot\delta_{y_i})=
\delta_{\bigoplus\limits_{i=1}^s\alpha_i\odot\delta_{y_i}}$ be an
arbitrary idempotent probability measure with everywhere finite
support from $\eta_{1,\ 2}(I(X))$. We have
\begin{center}
$H_{I^2}(I(\eta_{0,\ 1})(\mu),\
N)=H_{I^2}(\bigoplus\limits_{l=1}^n\lambda_l\odot\delta_{\delta_{x_l}},\
\delta_{\bigoplus\limits_{i=1}^s\alpha_i\odot\delta_{y_i}})=\bigoplus\limits_{l=1}^n
|\lambda_l|\odot H_I(\delta_{x_l},\
\bigoplus\limits_{i=1}^s\alpha_i\odot\delta_{y_i})=$\\
$=\bigoplus\limits_{l=1}^n |\lambda_l|\odot
\bigoplus\limits_{i=1}^s |\alpha_i|\odot \rho(x_l,\
y_i)=\bigoplus\limits_{i=1}^s |\alpha_i|\odot
\bigoplus\limits_{l=1}^n |\lambda_l|\odot \rho(x_l,\ y_i)=
\bigoplus\limits_{i=1}^s |\alpha_i|\odot
H_I(\bigoplus\limits_{l=1}^n \lambda_l\odot \delta_{x_l},\
\delta_{y_i})=$\\
$=\bigoplus\limits_{i=1}^s |\alpha_i|\odot H_I(\mu,\
\delta_{y_i})\geq \bigoplus\limits_{i=1}^s |\alpha_i|\odot \min
\{H_I(\mu,\ \delta_{y}):\ \delta_y\in \eta_{0,\ 1}(I(X))\}=$\\
$=\bigoplus\limits_{i=1}^s |\alpha_i|\odot H_I(\mu,\ \eta_{0,\
1}(I(X))\}\geq \bigoplus\limits_{i=1}^s |\alpha_i|\odot
\varepsilon\geq \varepsilon$.
\end{center}
Lemma 3 is proved.

Now Theorem 1 [3], Lemmas 1, 2 and 3 imply the main result of the
paper.

\textbf{Theorem 1.} \textsl{The functor $I$ is perfect
metrizable.}

Note that undelivered notions one can find in [5].

\end{document}